\documentclass[10pt]{article}
\usepackage{amsfonts}
\usepackage{mathrsfs}
\usepackage{amsthm}
\usepackage{amssymb}
\usepackage{amsmath}
\usepackage{enumerate}
\usepackage{graphicx}
\usepackage{subfigure}
\usepackage{float}
\usepackage{amsmath}
\numberwithin{equation}{section}

\theoremstyle{plain}
\newtheorem{theorem}{Theorem}[section]
\newtheorem{proposition}{Proposition}[section]

\newtheorem{lemma}[theorem]{Lemma}

\theoremstyle{definition}

\theoremstyle{remark}

\theoremstyle{Hypothesis}

\title{Least Squares Estimator for Vasicek Model Driven by Sub-fractional Brownian Processes from Discrete Observations }

\author { Cuiyun~Zhang,~\quad Jingjun Guo\footnote{Corresponding author (Jingjun Guo, guojj$@$Lzufe.edu.cn)},~\quad Aiqin Ma,~\quad Bo Peng\\
 School of Statistics,\\
Lanzhou University of Finance and Economics,\\
Lanzhou, Gansu 730020, PR China.\footnote{Research supported by the National Natural Science Foundation of China under Grant 71561017 and 71961013.}\\
 }

\begin{document}
\maketitle
\noindent\textbf{Abstract:}\ \ We study the parameter estimation problem of Vasicek Model driven by sub-fractional Brownian processes from discrete observations, and let $\{S_t^H, t\geq0\}$ denote a sub-fractional Brownian motion whose Hurst parameter $H\in(\frac{1}{2},1).$ The studies are as follows: firstly, two unknown parameters in the model are estimated by the least squares method. Secondly, the strong consistency and the asymptotic distribution of the estimators are studied respectively. Finally, our estimators are validated by numerical simulation.

\vskip 2mm
\noindent\textbf{\textit{$ \mathbf{Keywords}: $}}\ \ least squares method, Vasicek model, strong consistency, asymptotic distribution
\vskip 2mm

\section{Introduction}\label{sec-1}

\quad\quad The following Vasicek (1977) model driven by standard Brownian motion $\{B(t),~t\geq0\}$ has been extensively applied in various fields, such as economics, finance and environmental et al:
\begin{equation*}
dX_t=(\mu+\theta X_t)dt+\sigma dB_t^H,~~~t\geq 0,
\end{equation*}
where $\mu,~\theta$ are unknown parameters. The first term $(\mu+\theta X_t)$ is called the drift component, whose economic interpretation is that stochastic price fluctuations around the mean and price peaks are only temporarily, such as caused by power plant outages or capacity shortages.

Many extensions to this model have been made. For example, motivated by the phenomenon of long-range dependence found in data of telecommunication, economics and finance, the Brownian motion in the Vasicek model has been replaced by fractional Brownian motion (fBm). The fractional Vasicek model (fVm) was first used to describe the dynamics in volatility by Comte et al. (1998). Although fVm has many practical applications, little attention has been paid to its estimation and asymptotic theory in the literature. Xiao et al. (2019) developed the asymptotic theory for estimators of two parameters in the fVm. Tanaka et al. (2019) was concerned about the maximum likelihood method (MLE) of the drift parameters in the fVm from continuous observations.

Although this model driven by fBm has been applied in different areas, some more general fractional Gaussian processes, such as sub-fractional Brownian motion (sub-fBm), are still proposed. However, compared with the extensive studies of fBm, there are few systematic studies on statistical inference of other fractional Gaussian processes. The main reason for this phenomenon is the complexity of dependence structures fractional Gaussian processes which do not have stationary increments. Li et al.(2018) tackled the least squares estimators (LSE) and discussed the consistency and asymptotic distributions of the two estimators in the Vasicek model driven by sub-fBm based on the continuous observations. Xiao et al. (2018) considered the parameter estimation for the continuously observed Vasicek model with sub-fBm. Furthermore, the strong consistency results as well as the asymptotic distributions of these estimators are obtained in both the non-ergodic case and the null recurrent case.

From a practical point of view, it is more realistic and interesting to consider parameter estimation based on discrete observations in statistical inference, and the asymptotic theory of parameter estimation for stochastic processes is also well developed. Shen et al. (2020) considered the problem of parameter estimation for Vasicek model driven by small fractional L\'{e}vy noise based on discrete high-frequency observations at regularly spaced time points. For the general case and the null recurrent case, the consistency as well as the asymptotic behavior of LSE of two unknown parameters have been established.

Motivated by the aforementioned works, in this article, we study the LSE for Vasicek model:
\begin{equation}\label{eq-1}
dX_t=(\mu+\theta X_t)dt+\sigma dS_t^H,~~t\geq 0,~~X_0=x_0,
\end{equation}
where $S_t^H$ is a sub-fBm with Hurst index $H\in(\frac{1}{2},1),$ $x_0$ is a fix value. In almost all empirically relevant cases, the parameters $\mu$ and $\theta$ in the drift component of model (\ref{eq-1}) are unknown and the real value of these two parameters are $\theta_0$ and $\mu_0.$ We assume to observe $\{X_t,t\geq0\}$ at $n$ regular time intervals $\{t_i=\frac{i}{n},i=1,2,\cdots, n\}$, so an important problem is to estimate parameters $\theta$ and $\mu$ according to $\{X_t,t\geq0\}.$

The rest of the paper is organized as follows. In section 2, we introduce the detailed information about sub-fBm in preparation for our proof and describe the LSE of Vasicek model driven by sub-fBm from discrete observations. The strong consistency of LSE for our model are given in Section 3. Section 4 is devoted to the asymptotic distribution of LSE for Vasicek model. In section 5, our estimations are validated by numerical simulations. The true values of the parameters are given and then they are used to simulate the Vasicek model driven by sub-fBm. With these simulated values we compute our estimators and compare them with the true parameters. Numerical results show that our estimators converges to the true parameters.

\section{ Preliminaries }\label{sec-2}

\quad\quad In this section, we describe some basic facts of sub-fBm and the LSE of Vasicek model driven by sub-fBm from discrete observations. More complete introductions to this subjects, see Mendy (2013), Nourdin et al. (2017), Tudor (2007) and the references therein.

The sub-fBm arises from occupation time fluctuations of branching particle systems with Poisson initial condition. As we all know, the sub-fBm has properties similar to fBm such as self-similarity, long-range dependence and H\"{o}lder continuous paths. However, compared with fBm, sub-fBm has non stationary increments. The increments over non overlapping intervals are more weakly correlated and their covariance decays polynomially at a higher rate. For this reason, it is called sub-fBm in Bojdecki et al.(2004). It is worth emphasizing that the properties mentioned here make the sub-fBm a possible candidate for models involving long-range dependence, self-similarity and non-stationary.

The sub-fBm $S_t^H$ is a mean zero Gaussian process with $S_0^H=0$ and the covariance
\begin{align*}
C_H(s,t)=\mathbf{E}(S_t^HS_s^H)=s^{2H}+t^{2H}-\frac{1}{2}\{\mid s-t\mid^{2H}+(s+t)^{2H}\},
\end{align*}
where $s,~t\geq 0.$ When $H=\frac{1}{2},$ $S_t^H$ coincides with the standard Brownian motion.  Actually, $S_t^H$ is neither a semimartingale nor a Markov process unless $H=\frac{1}{2}.$

For all $s\leq t,$ there is
\begin{equation}\label{eq-2}
\mathbf{E}(\mid S_t^H-S_s^H\mid^2)=-2^{2H-1}(t^{2H}+s^{2H})+(t+s)^{2H}-(t-s)^{2H}.
\end{equation}
The increments of sub-fBm satisfy the following inequalities
\begin{equation}\label{eq-3}
[(2-2^{2H-1})\wedge 1](t-s)^{2H}\leq \mathbf{E}(\mid S_t^H-S_s^H\mid^2)\leq [(2-2^{2H-1})\vee 1](t-s)^{2H}.
\end{equation}
Moreover, for $u\leq v\leq s\leq t$ the covariance of increments of sub-fBm over non-overlapping intervals can be written as
\begin{align*}
\mathbf{E}((S_t^H-S_s^H)(S_v^H-S_u^H))=\frac{1}{2} [(t+u)^{2H}+(t-u)^{2H}+(s+v)^{2H}+(s-v)^{2H}\\
-(t+v)^{2H}-(t-v)^{2H}-(s+u)^{2H}-(s-u)^{2H}].
\end{align*}

Fixed a time interval $[0,T],$ We denote by $\mathcal{H}_{S_t^H}$ canonical Hilbert space associated to the sub-fBm $S_t^H$. That is, $\mathcal{H}_{S_t^H}$ is the closure of the linear span $\varepsilon$ generated by the indicator function with respect to the scalar product
\begin{align*}
\langle\mathbf{I}_{[0,t]},\mathbf{I}_{[0,s]}\rangle_{\mathcal{H}_{S_t^H}} =C_H(s,t).
\end{align*}
The covariance of sub-fBm also can be written as
\begin{align*}
C_H(s,t)=\mathbf{E}(S_t^HS_s^H)=\int_0^t\int_0^s\phi_H(u,v)dudv,
\end{align*}
where $\phi_H(u,v)=H(2H-1)[\mid t-s\mid^{2H-2}-(t+s)^{2H-2}]$ and $\frac{1}{2}<H<1.$

For $H>\frac{1}{2},$ we have $L^{\frac{1}{H}}([0,T])\subset \mathcal{H}_{S_t^H}$ and for any pair step function $\varphi,\psi\in L^{\frac{1}{H}}([0,T]):$
\begin{equation}\label{eq-4}
\langle\varphi,\psi\rangle_\mathcal{H}=\alpha_H\int_0^T\int_0^T\varphi_s\psi_t\phi_H(s,t)dsdt,~~~\frac{1}{2}<H<1.
\end{equation}

Next, let's consider the Vasicek model driven by sub-fBm, which takes the sub-fBm as the governing force of the state variable instead of the usual Brownian motion.

For the stochastic differential equation (\ref{eq-1}), we discuss the LSE of the two parameters.

LSE's motivation is the following illuminating argument, minimizing contrast function of $\mu$ and $\theta$ respectively,
\begin{align*}
\rho_{n,\sigma}(\theta,\mu)=\sum_{i=1}^{n}\mid X_{t_i}-X_{t_{i-1}}-(\mu+\theta X_{t_{i-1}})\cdot\bigtriangleup t_{i-1}\mid^2,
\end{align*}
where $\bigtriangleup t_{i}=t_i-t_{i-1}=\frac{1}{n},i=1,2,\cdots, n.$

Taking the partial derivatives of $\theta$ and $\mu$ respectively, we get
\begin{align*}
&\frac{\partial \rho_{n,\sigma}(\theta,\mu)}{\partial \mu}=\sum\limits_{i=1}^{n}( X_{t_i}-X_{t_{i-1}}-(\mu+\theta X_{t_{i-1}})\cdot\bigtriangleup t_{i-1})=0,\\
&\frac{\partial \rho_{n,\sigma}(\theta,\mu)}{\partial \theta}=\sum\limits_{i=1}^{n}( X_{t_i}-X_{t_{i-1}}-(\mu+\theta X_{t_{i-1}})\cdot\bigtriangleup t_{i-1})X_{t_{i-1}}=0.
\end{align*}

To solve the above equation, we have
\begin{align}\label{eq-5}
\hat{\theta}=\frac{\sum\limits_{i=1}^{n}(X_{t_i}-X_{t_{i-1}})X_{t_{i-1}}-\frac{1}{n}\sum\limits_{i=1}^{n}X_{t_{i-1}}\sum\limits_{i=1}^{n}(X_{t_i}-X_{t_{i-1}})}
{\frac{1}{n}\sum\limits_{i=1}^{n}X_{t_{i-1}}^2-\frac{1}{n^2}(\sum\limits_{i=1}^{n}X_{t_{i-1}})^2},
\end{align}
\begin{align}\label{eq-6}
\hat{\mu}=\frac{\sum\limits_{i=1}^{n}(X_{t_i}-X_{t_{i-1}})\sum\limits_{i=1}^{n}X_{t_{i-1}}^2-\sum\limits_{i=1}^{n}X_{t_{i-1}}\sum\limits_{i=1}^{n}(X_{t_i}-X_{t_{i-1}})X_{t_{i-1}}}
{\sum\limits_{i=1}^{n}X_{t_{i-1}}^2-\frac{1}{n}(\sum\limits_{i=1}^{n}X_{t_{i-1}})^2}.
\end{align}

\section{The Consistency of the Least Squares Estimator}\label{sec-3}
\quad\quad In this section, our main purpose is to clarify and prove the Theorem 3.1, which gives the consistency of the estimators given by equations (\ref{eq-5}) and (\ref{eq-6}).

Let's consider the following solution of the stochastic differential equation (\ref{eq-1}):
\begin{align}\label{eq-7}
X_t=x_0e^{\theta t}+\frac{\mu}{\theta}(e^{\theta t}-1)+\sigma\int_0^te^{\theta (t-s)}dS_s^H, ~~~~~~t\in[0,1].
\end{align}
More specifically, the numerical approximation of the model (\ref{eq-1}) can be expressed as Eulerian model (Ait-Sahalia (2002)):
\begin{align}\label{eq-8}
X_{t_i}=X_{t_{i-1}}+(\mu+\theta X_{t_{i-1}})\bigtriangleup t_i+\sigma(S_{t_i}^H-S_{t_{i-1}}^H).
\end{align}

Hence, substituting (\ref{eq-8}) into (\ref{eq-5}) and (\ref{eq-6}) respectively, we get:
\begin{align}
\hat{\theta}=\theta_0+\sigma\frac{\sum\limits_{i=1}^{n}X_{t_{i-1}}(S_{t_i}^H-S_{t_{i-1}}^H)-\frac{1}{n}\sum\limits_{i=1}^{n}X_{t_{i-1}}S_n^H}
{\frac{1}{n}\sum\limits_{i=1}^{n}X_{t_{i-1}}^2-\frac{1}{n^2}(\sum\limits_{i=1}^{n}X_{t_{i-1}})^2},
\end{align}
\begin{align}
\hat{\mu}=\mu_0+\sigma\frac{\sum\limits_{i=1}^{n}X_{t_{i-1}}^2S_n^H-\sum\limits_{i=1}^{n}X_{t_{i-1}}\sum\limits_{i=1}^{n}X_{t_{i-1}}(S_{t_i}^H-S_{t_{i-1}}^H)}
{\sum\limits_{i=1}^{n}X_{t_{i-1}}^2-\frac{1}{n}(\sum\limits_{i=1}^{n}X_{t_{i-1}})^2},
\end{align}
where $\theta_0,~\mu_0$ are the true values of parameters $\mu$ and $\theta$ respectively.

Next, we will state our main results.
\begin{theorem}\label{Th-3-1} For $H\in(\frac{1}{2},1)$, we have

(1) $\hat{\theta}~\underrightarrow{\text{a.s}}~\theta_0$, as $n\rightarrow\infty$ and $\sigma\rightarrow0;$

(2) $\hat{\mu}~\underrightarrow{\text{a.s}}~\mu_0,$ as $n\rightarrow\infty$ and $\sigma\rightarrow0.$
\end{theorem}

In order to simplify the proof of Theorem 3.1, we firstly give the following lemmas and propositions.

For simplicity, we assume that
\begin{align}
X_t^0=\frac{\mu}{\theta}(e^{\theta t}-1)+x_0e^{\theta t},~~t\in[0,1].
\end{align}

\begin{lemma}\cite{21}\label{lem-3-2} For any $0<u_2\leq u_1\leq v_1,~0<u_2\leq v_2\leq v_1$ and $u_1-u_2=v_1-v_2,$ there exists a constant $C$ depend on $\theta$ and $H$ such that
\begin{align*}
|\int_{u_2}^{u_1}\int_{v_2}^{v_1}e^{-\theta(s+t)}|s-t|^{2H-2}dsdt|\leq C|e^{-\theta(u_1+v_1)}-e^{-\theta(u_2+v_2)}||v_1-u_2|^{2H-1},   ~~~\theta\neq0.
\end{align*}
\end{lemma}
\begin{lemma}\label{lem-3-3} For $\theta < 0,$ we have
\begin{align*}
&\mathbf{E}(\int_{t_{i-1}}^{t_i}e^{\theta(t_i-s)}dS_s^H)^2\leq C|e^{\frac{-2\theta }{n}}-1|n^{1-2H},\\
&\mathbf{E}(\int_{0}^{t_{i-1}}e^{\theta(t_{i-1}-s)}dS_s^H)^2\leq C|e^{-2\theta t_{i-1}-1 }||t_{i-1}|^{2H-1},
\end{align*}
where the $C$ depend on $H$ and $\theta.$
\end{lemma}
\noindent\textbf{Proof}:From (\ref{eq-4}), we calculate the following formula directly:
\begin{align*}
&\mathbf{E}(\int_{t_{i-1}}^{t_i}e^{\theta(t_i-s)}dS_s^H)^2=H(2H-1)\int_{t_{i-1}}^{t_i}\int_{t_{i-1}}^{t_i}e^{\theta(2t_i-u-v)}(|u-v|^{2H-2}-(u+v)^{2H-2})dudv\\
&=H(2H-1)\int_{t_{i-1}}^{t_i}\int_{t_{i-1}}^{t_i}e^{\theta(2t_i-u-v)}|u-v|^{2H-2}dudv\\
&-H(2H-1)\int_{t_{i-1}}^{t_i}\int_{t_{i-1}}^{t_i}e^{\theta(2t_i-u-v)}(u+v)^{2H-2}dudv\\
&\leq H(2H-1)\int_{t_{i-1}}^{t_i}\int_{t_{i-1}}^{t_i}e^{\theta(2t_i-u-v)}|u-v|^{2H-2}dudv
\end{align*}
According to lemma 3.2, the right side of the above inequality satisfies the following inequality
\begin{align*}
&H(2H-1)\int_{t_{i-1}}^{t_i}\int_{t_{i-1}}^{t_i}e^{\theta(2t_i-u-v)}|u-v|^{2H-2}dudv\\
&\leq H(2H-1)e^{2\theta}\int_{t_{i-1}}^{t_i}\int_{t_{i-1}}^{t_i}e^{-\theta(u+v)}|u-v|^{2H-2}dudv\\
&\leq C|e^{\frac{-2\theta }{n}}-1|n^{1-2H}
\end{align*}
So, we obtain
\begin{align*}
\mathbf{E}(\int_{t_{i-1}}^{t_i}e^{\theta(t_i-s)}dS_s^H)^2\leq C|e^{\frac{-2\theta }{n}}-1|n^{1-2H}.
\end{align*}
As the same proof method as the above inequality, we can get
\begin{align*}
&\mathbf{E}(\int_0^{t_{i-1}}e^{\theta({t_{i-1}}-s)}dS_s^H)^2=\\
&H(2H-1)\int_0^{t_{i-1}}\int_0^{t_{i-1}}e^{\theta(2t_{i-1}-u-v)}(|u-v|^{2H-2}-(u+v)^{2H-2})dudv\\
&\leq H(2H-1)\int_0^{t_{i-1}}\int_0^{t_{i-1}}e^{\theta(2t_{i-1}-u-v)}|u-v|^{2H-2}dudv\\
&\leq H(2H-1)e^{2\theta}\int_0^{t_{i-1}}\int_0^{t_{i-1}}e^{-\theta(u+v)}|u-v|^{2H-2}dudv\\
&\leq C|e^{-2\theta t_{i-1})}-1||t_{i-1}|^{2H-1}.\qed
\end{align*}

\begin{proposition}\label{pro-3-1} As $\sigma\rightarrow0,$ we have
\begin{align*}
\sup\limits_{t\in[0,1]}\mid(X_t)^2-(X_t^0)^2\mid\rightarrow0,
\end{align*}
where $X_t^0$ is the equation $(3.1).$
\end{proposition}
\noindent\textbf{Proof}: For $(1.1),$ we can rewrite as
\begin{align*}
X_t=X_0+\int_0^t(\mu+\theta X_s)ds+\sigma S_t^H,~~~t\in[0,1].
\end{align*}
Then, from Equation $(3.5),$ we have
\begin{align*}
X_t-X_t^0=\int_0^t\theta(X_s-X_s^0)ds+\sigma S_t^H,~~~t\in[0,1].
\end{align*}

On the one hand, by Cauchy-Schwarz inequality, we get
\begin{align*}
\mid X_t-X_t^0\mid^2&\leq 2\mid\int_0^t\theta(X_s-X_s^0)ds\mid^2+2\sigma^2 \mid S_t^H\mid^2\\
&\leq 2t\int_0^t\mid\theta(X_s-X_s^0)\mid^2ds+2\sigma^2 \mid S_t^H\mid^2\\
&\leq 2\theta^2 t\int_0^t\mid(X_s-X_s^0)\mid^2ds+2\sigma^2 \mid S_t^H\mid^2.
\end{align*}
And by Gronwall's inequality, we get the following inequality
\begin{align*}
\mid X_t-X_t^0\mid^2\leq2\sigma^2e^{2\theta^2 t^2}\mid S_t^H\mid^2.
\end{align*}
Then,
\begin{align*}
\mid X_t-X_t^0\mid^2\leq2\sigma^2e^{2\theta^2 t^2}\sup\limits_{0\leq s\leq t}\mid S_s^H\mid^2.
\end{align*}
Therefore, we find
\begin{align*}
\sup\limits_{0\leq t\leq 1}\mid X_t-X_t^0\mid\leq\sqrt{2}\sigma e^{\theta^2}\sup\limits_{0\leq t\leq 1}\mid S_t^H\mid^2.
\end{align*}
So in summary, we can see
\begin{align}
\sup\limits_{0\leq t\leq 1}\mid X_t-X_t^0\mid\rightarrow 0,  ~~~\sigma\rightarrow 0.
\end{align}

On the other hand, according to the same method, we have
\begin{align*}
\mid X_t\mid^2&=\mid X_0+\int_0^t(\mu+\theta X_s)ds+\sigma S_t^H\mid^2\\
&\leq2(\mid X_0\mid+\sigma\mid S_t^H\mid+\mid\mu t\mid)^2+2\mid\int_0^t\theta X_sds\mid^2\\
&\leq2(\mid X_0\mid+\sigma\sup\limits_{0\leq t\leq 1}\mid S_t^H\mid+\mid\mu\mid)^2+2\theta^2t\int_0^t\mid X_s\mid^2ds.
\end{align*}
By Gronwall's inequality, we have
\begin{align*}
\mid X_t\mid\leq\sqrt{2}(\mid x_0\mid+\sigma\sup\limits_{0\leq t\leq 1}\mid S_t^H\mid+\mid\mu\mid)e^{\theta^2t^2}<\infty.
\end{align*}
Thus
\begin{align*}
&\sup\limits_{0\leq t\leq 1}\mid(X_t)^2-(X_t^0)^2\mid\\
&\leq(\sup\limits_{0\leq t\leq 1}\mid X_t\mid+\sup\limits_{0\leq t\leq 1}\mid X_t^0\mid)(\sup\limits_{0\leq t\leq 1}\mid X_t-X_t^0)\rightarrow 0,~~~~~~~~~~\sigma\rightarrow 0.\qed
\end{align*}

\begin{proposition}\label{pro-3-5} As $\sigma\rightarrow 0$ and $n\rightarrow \infty,$ then we have
\begin{align*}
\frac{1}{n}\sum_{i=1}^{n}X_{t_{i-1}}^2\rightarrow \int_0^1(X_t^0)^2dt, \\
\frac{1}{n}\sum_{i=1}^{n}X_{t_{i-1}}\rightarrow \int_0^1(X_t^0)dt.
\end{align*}
\end{proposition}
\noindent\textbf{Proof}: Owing to the following equation:
\begin{align*}
\frac{1}{n}\sum_{i=1}^{n}X_{t_{i-1}}^2&=\sum_{i=1}^n\int_{t_{i-1}}^{t_i}X_{t_{i-1}}^2dt
=\sum_{i=1}^n\int_{t_{i-1}}^{t_i}X_{\frac{[nt]}{n}}^2dt\\
&=\int_0^1(X_{\frac{[nt]}{n}})^2dt,
\end{align*}
where $t_i-t_{i-1}=\frac{1}{n},$ $[nt]$ denotes the integer part of $nt$. Then
\begin{align*}
\int_0^1(X_{\frac{[nt]}{n}})^2dt\rightarrow\int_0^1(X_t)^2dt,~~~~~~as ~~~n\rightarrow\infty.
\end{align*}
By Proposition 3.1, there is
\begin{align*}
&\frac{1}{n}\sum_{i=1}^{n}X_{t_{i-1}}^2-\int_0^1(X_t^0)^2dt=\int_0^1(X_{\frac{[nt]}{n}})^2dt-\int_0^1(X_t^0)^2dt\\
&\leq\sup\limits_{0\leq t\leq 1}\mid (X_{\frac{[nt]}{n}})^2-X_{\frac{[nt]}{n}}^0)^2\mid+\sup\limits_{0\leq t\leq 1}\mid (X_{\frac{[nt]}{n}}^0)^2-(X_t^0)^2\mid\\
&\rightarrow 0,~~~~~~~~~~~~n\rightarrow\infty.
\end{align*}
In the same way, we can obtain
\begin{align*}
 \frac{1}{n}\sum_{i=1}^{n}X_{t_{i-1}}\rightarrow \int_0^1(X_t^0)dt, ~~~~~~~~~~~n\rightarrow\infty.\qed
\end{align*}

\begin{lemma}\label{lem-3-6} For $n\rightarrow\infty,~\sigma\rightarrow0,$ then
\begin{align*}
\sum_{i=1}^{n}(X_{t_{i-1}}(S_{t_i}^H-S_{t_{i-1}}^H))\rightarrow x_0\int_0^1e^{\theta t}dS_t^H+\frac{\mu}{\theta}\int_0^1(e^{\theta t}-1)dS_t^H<\infty.
\end{align*}
\end{lemma}
\noindent\textbf{Proof}: According to the solution of (\ref{eq-1}), we have
\begin{align*}
X_{t_{i-1}}=x_0e^{\theta t_{i-1}}+\frac{\mu}{\theta}(e^{\theta t_{i-1}}-1)+\sigma\int_0^{t_{i-1}}e^{\theta (t_{i-1}-s)}dS_s^H,
\end{align*}
where $t_i-t_{i-1}=\frac{1}{n},~~i=1,2,\cdots, n.$ Hence
\begin{align*}
&\sum_{i=1}^{n}X_{t_{i-1}}(S_{t_i}^H-S_{t_{i-1}}^H)=\sum_{i=1}^{n}(x_0e^{\theta t_{i-1}}+\frac{\mu}{\theta}(e^{\theta t_{i-1}}-1)+\sigma\int_0^{t_{i-1}}e^{\theta (t_{i-1}-s)}dS_s^H)(S_{t_i}^H-S_{t_{i-1}}^H)\\
&=I_1+I_2+I_3.
\end{align*}
For $\theta<0,$ we can rewrite as
\begin{align*}
I_1=\sum_{i=1}^{n}x_0e^{\theta t_{i-1}}(S_{t_i}^H-S_{t_{i-1}}^H)=x_0\int_0^1e^{\theta t}dS_t^H.
\end{align*}

\begin{align*}
I_2=\sum_{i=1}^{n}\frac{\mu}{\theta}(e^{\theta t_{i-1}}-1)(S_{t_i}^H-S_{t_{i-1}}^H)=\frac{\mu}{\theta}\int_0^1(e^{\theta t}-1)dS_t^H.
\end{align*}

For $I_3,$ by the Markov inequality, there exists any $\delta>0$
\begin{align*}
&P(\mid\sigma\sum\limits_{i=1}^n\int_0^{t_{i-1}}e^{\theta (t_{i-1}-s)}dS_s^H(S_{t_i}^H-S_{t_{i-1}}^H)\mid>\delta)\\
&\leq \delta^{-1}\sigma\sum\limits_{i=1}^n(\mathbf{E}(\int_0^{t_{i-1}}e^{\theta (t_{i-1}-s)}dS_s^H)^2)^{\frac{1}{2}}(\mathbf{E}(S_{t_i}^H-S_{t_{i-1}}^H)^2)^{\frac{1}{2}}\\
&\leq C\delta^{-1}\sigma\sum\limits_{i=1}^n|e^{-2\theta t_{i-1}}-1|^{\frac{1}{2}}|t_{i-1}|^{H-\frac{1}{2}}(|t_i-t_{i-1}|^{2H})^{\frac{1}{2}}\\
&\leq C\delta^{-1}\sigma\sum\limits_{i=1}^n|e^{-2\theta t_{i-1}}-1|^{\frac{1}{2}}n^{-H}\rightarrow ~0,   ~~(n\rightarrow\infty,~\sigma\rightarrow0),
\end{align*}
where $C$ is a constant depend on $H$ and $\theta,$ so $I_3\rightarrow 0$ as $n\rightarrow\infty,~\sigma\rightarrow0.$ \qed

\noindent\textbf{Proof of Theorem 3.1}: Combined propositions $3.1~3.2$ with Lemma $3.4,$ when $n\rightarrow\infty,~~\sigma\rightarrow0,$ we have
\begin{align*}
&\sigma(\sum_{i=1}^{n}X_{t_i}(S_{t_i}^H-S_{t_{i-1}}^H)-\frac{1}{n}\sum_{i=1}^{n}X_{t_{i-1}}S_n^H)\rightarrow0,\\
&\frac{1}{n}\sum_{i=1}^{n}X_{t_{i-1}}^2-\frac{1}{n^2}(\sum_{i=1}^{n}X_{t_{i-1}})^2\rightarrow\int_0^1(X_t^0)^2dt-(\int_0^1X_t^0dt)^2.
\end{align*}
So, we immediately come to the conclusion:
when $n~\rightarrow ~\infty$ and $\sigma\rightarrow0$, then $\hat{\theta}~\underrightarrow{\text{a.s}}~\theta_0.$
Moreover,
\begin{align*}
\hat{\mu}=\mu_0+\sigma\frac{\sum\limits_{i=1}^{n}X_{t_{i-1}}^2S_n^H-\sum\limits_{i=1}^{n}X_{t_{i-1}}\sum\limits_{i=1}^{n}X_{t_{i-1}}(S_{t_i}^H-S_{t_{i-1}}^H)}
{\sum\limits_{i=1}^{n}X_{t_{i-1}}^2-\frac{1}{n}(\sum\limits_{i=1}^{n}X_{t_{i-1}})^2}
\end{align*}
\begin{align*}
=\mu_0+\sigma\frac{\frac{1}{n}\sum\limits_{i=1}^{n}X_{t_{i-1}}^2S_n^H-\frac{1}{n}\sum\limits_{i=1}^{n}X_{t_{i-1}}\sum\limits_{i=1}^{n}X_{t_{i-1}}(S_{t_i}^H-S_{t_{i-1}}^H)}
{\frac{1}{n}\sum\limits_{i=1}^{n}X_{t_{i-1}}^2-(\frac{1}{n}\sum\limits_{i=1}^{n}X_{t_{i-1}})^2}.
\end{align*}

Similarity, we can get $\hat{\mu}~\underrightarrow{\text{a.s}}~\mu_0,$ as $n\rightarrow\infty$ and $\sigma\rightarrow0.$ \qed

\section{Asymptotic Distributions of the LSE}\label{sec-4}
\quad\quad According to  Es-Sebaiy (2013), Wang et al. (2017) and the solution of equation (\ref{eq-1}), the process is observed at equidistant discrete times $\{t_i=\frac{i}{n},i=1,2,\cdots, n\},$ so we can obtain
\begin{align*}
X_{t_i}=X_{t_{i-1}}e^{\frac{\theta_0 }{n}}+\frac{\mu}{\theta_0}(e^\frac{\theta_0 }{n}-1)+\sigma\int_{t_{i-1}}^{t_i}e^{\theta_0 (t_i-s)}dS_s^H,
\end{align*}
where $t_i-t_{i-1}=\frac{1}{n}.$

Then the estimated value of the parameter $\theta, \mu$  can be rewritten as
\begin{align}
\hat{\theta}=\frac{e^{\frac{\theta_0}{n}}-1}{n^{-1}}+\sigma\frac{\sum\limits_{i=1}^{n}X_{t_{i-1}}\int_{t_{i-1}}^{t_i}e^{\theta_0(t_i-s)}dS_s^H-\frac{1}{n}\sum\limits_{i=1}^{n}X_{t_{i-1}}
\sum\limits_{i=1}^{n}\int_{t_{i-1}}^{t_i}e^{\theta_0(t_i-s)}dS_s^H}
{\frac{1}{n}\sum\limits_{i=1}^{n}X_{t_{i-1}}^2-\frac{1}{n^2}(\sum\limits_{i=1}^{n}X_{t_{i-1}})^2},
\end{align}
\begin{align}
\hat{\mu}=\frac{\frac{\mu_0}{\theta_0}(e^{\frac{\theta_0}{n}}-1)}{n^{-1}}+\sigma\frac{\sum\limits_{i=1}^{n}X_{t_{i-1}}^2\sum\limits_{i=1}^{n}\int_{t_{i-1}}^{t_i}e^{\theta_0(t_i-s)}dS_s^H
-\sum\limits_{i=1}^{n}X_{t_{i-1}}\sum\limits_{i=1}^{n}X_{t_{i-1}}\int_{t_{i-1}}^{t_i}e^{\theta_0(t_i-s)}dS_s^H}
{\sum\limits_{i=1}^{n}X_{t_{i-1}}^2-\frac{1}{n}(\sum\limits_{i=1}^{n}X_{t_{i-1}})^2}.
\end{align}

\begin{lemma}\label{lem-4-1} As $n\rightarrow\infty,$ then
\begin{align*}
\sum\limits_{i=1}^{n}X_{t_{i-1}}\int_{t_{i-1}}^{t_i}e^{\theta_0(t_i-s)}dS_s^H\rightarrow\int_0^1X_s^0dS_s^H.
\end{align*}
\end{lemma}

\noindent\textbf{Proof}: Using the same methods as the proof of lemma 3.4, we have
\begin{align*}
\sum\limits_{i=1}^{n}X_{t_{i-1}}\int_{t_{i-1}}^{t_i}e^{\theta_0(t_i-s)}dS_s^H&=\int_0^1\sum\limits_{i=1}^{n}X_{t_{i-1}}e^{\theta_0(t_i-s)}\mathbf{I}_{[t_{i-1},t_i]}dS_s^H\\
&=\int_0^1e^{\theta_0(\frac{[ns]+1}{n}-s)}X_{\frac{[ns]}{n}}dS_s^H.
\end{align*}

On the other hand,
\begin{align*}
e^{\theta_0(\frac{[ns]+1}{n}-s)}X_{\frac{[ns]}{n}}\rightarrow X_s,~~~n\rightarrow\infty.
\end{align*}
and according to equation $(3.6),$ $X_s\rightarrow X_s^0,~~~~~~~~~~n\rightarrow\infty.$

Therefore, we have
\begin{align*}
\sum\limits_{i=1}^{n}X_{t_{i-1}}\int_{t_{i-1}}^{t_i}e^{\theta_0(t_i-s)}dS_s^H\rightarrow\int_0^1X_s^0dS_s^H,~~~~~~~~~~~~n\rightarrow\infty.\qed
\end{align*}

\begin{theorem}\label{Th-4-2} As $n\rightarrow\infty,~~\sigma\rightarrow 0$ and $n\sigma\rightarrow\infty,$ we obtain
\begin{align*}
\sigma^{-1}(\hat{\theta}-\theta_0)\rightarrow \frac{\int_0^1X_s^0dS_s^H-\int_0^1X_s^0ds\int_0^1dS_s^H}{\int_0^1(X_s^0)^2ds-(\int_0^1X_s^0ds)^2},
\end{align*}
\begin{align*}
\sigma^{-1}(\hat{\mu}-\mu_0)\rightarrow \frac{\int_0^1(X_s^0)^2ds\int_0^1dS_s^H-\int_0^1X_s^0ds\int_0^1X_s^0dS_s^H}{\int_0^1(X_s^0)^2ds-(\int_0^1X_s^0ds)^2}.
\end{align*}
\end{theorem}

\noindent\textbf{Proof}: According to $(4.1)$
\begin{align*}
\sigma^{-1}(\hat{\theta}-\theta_0)=\sigma^{-1}(\frac{e^{\frac{\theta_0}{n}}-1}{n^{-1}}-\theta_0)+\frac{\sum\limits_{i=1}^{n}X_{t_i}\int_{t_{i-1}}^{t_i}e^{\theta_0(t_i-s)}dS_s^H-\frac{1}{n}\sum\limits_{i=1}^{n}X_{t_{i-1}}
\sum\limits_{i=1}^{n}\int_{t_{i-1}}^{t_i}e^{\theta_0(t_i-s)}dS_s^H}
{\frac{1}{n}\sum\limits_{i=1}^{n}X_{t_{i-1}}^2-\frac{1}{n^2}(\sum\limits_{i=1}^{n}X_{t_{i-1}})^2},
\end{align*}
Obviously, $\sigma^{-1}(\frac{e^{\frac{\theta_0}{n}}-1}{n^{-1}}-\theta_0)\rightarrow 0,$ if $n\rightarrow\infty$ and $n\sigma\rightarrow\infty.$
Moreover,
\begin{align*}
&\sum\limits_{i=1}^{n}\int_{t_{i-1}}^{t_i}e^{\theta_0(t_i-s)}dS_s^H=\int_0^1\sum\limits_{i=1}^{n}e^{\theta_0(t_i-s)}dS_s^H\\
&=\int_0^1e^{\theta_0(\frac{[ns]}{n}-s)}dS_s^H=\int_0^1dS_s^H.
\end{align*}
Combining Lemma 4.1 with Proposition 3.2, as $n\rightarrow\infty,~~\sigma\rightarrow 0$ and $n\sigma\rightarrow\infty,$ we obtain
\begin{align*}
\sigma^{-1}(\hat{\theta}-\theta_0)\rightarrow\frac{\int_0^1X_s^0dS_s^H-\int_0^1X_s^0ds\int_0^1dS_s^H}{\int_0^1(X_s^0)^2ds-(\int_0^1X_s^0ds)^2}
\end{align*}
Further, we calculate the following equation
\begin{align*}
\sigma^{-1}(\hat{\mu}-\mu_0)=&\sigma^{-1}(\frac{\frac{\mu_0}{\theta_0}(e^{\frac{\theta_0}{n}}-1)}{n^{-1}}-\mu_0)\\
&+\frac{\sum\limits_{i=1}^{n}X_{t_{i-1}}^2\sum\limits_{i=1}^{n}\int_{t_{i-1}}^{t_i}e^{\theta_0(t_i-s)}dS_s^H
-\sum\limits_{i=1}^{n}X_{t_{i-1}}\sum\limits_{i=1}^{n}X_{t_{i-1}}\int_{t_{i-1}}^{t_i}e^{\theta_0(t_i-s)}dS_s^H}
{\sum\limits_{i=1}^{n}X_{t_{i-1}}^2-\frac{1}{n}(\sum\limits_{i=1}^{n}X_{t_{i-1}})^2}\\
=&\sigma^{-1}(\frac{\frac{\mu_0}{\theta_0}(e^{\frac{\theta_0}{n}}-1)}{n^{-1}}-\mu_0)\\
&+\frac{\frac{1}{n}\sum\limits_{i=1}^{n}X_{t_{i-1}}^2\sum\limits_{i=1}^{n}\int_{t_{i-1}}^{t_i}e^{\theta_0(t_i-s)}dS_s^H
-\frac{1}{n}\sum\limits_{i=1}^{n}X_{t_{i-1}}\sum\limits_{i=1}^{n}X_{t_{i-1}}\int_{t_{i-1}}^{t_i}e^{\theta_0(t_i-s)}dS_s^H}
{\frac{1}{n}\sum\limits_{i=1}^{n}X_{t_{i-1}}^2-(\frac{1}{n}\sum\limits_{i=1}^{n}X_{t_{i-1}})^2}.
\end{align*}

In the same way, combining Lemma 4.1 with Proposition 3.2, as $n\rightarrow\infty,~~\sigma\rightarrow 0$ and $n\sigma\rightarrow\infty,$ we get
\begin{align*}
\sigma^{-1}(\hat{\mu}-\mu_0)\rightarrow\frac{\int_0^1(X_s^0)^2ds\int_0^1dS_s^H-\int_0^1X_s^0ds\int_0^1X_s^0dS_s^H}{\int_0^1(X_s^0)^2ds-(\int_0^1X_s^0ds)^2}.\qed
\end{align*}

\section{Simulation}\label{sec-5}
\quad\quad In this section, we use Monte Carlo simulation to prove the unbiasedness and effectiveness of estimators $\theta$ and $\mu.$

First of all, we use R to simulate stochastic differential equation $(1.1)$, as show in Figure 1. The unbiasedness and effectiveness of the estimators are validated respectively in Table 1 and Table 2.

Next, to calculate the mean value and standard deviation of the estimators $\theta$ and $\mu$ (as shown in Table 1 and Table 2, respectively), we let $\sigma=0.4,$ $x_0=0$ take different values of $H,$ $\theta$ and use R to generate 500 samples according to the equation $(3.3)$ and $(3.4)$. We can see that the mean value almost converges to the real value and the standard deviation is relatively small, which shows that this estimators are relatively accurate.
\begin{figure}[H]
\centering
\includegraphics[scale=1]{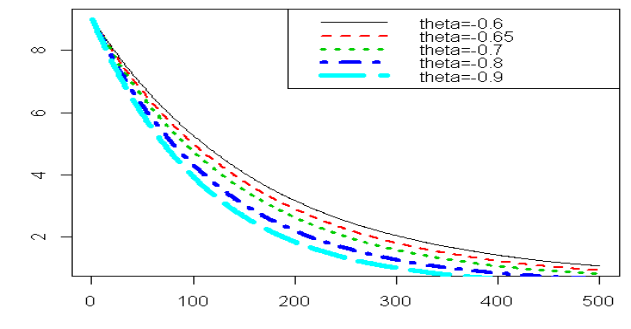}
\caption{The figure of  $X_t,$ where $x_0=9,~\mu=0.4,~\sigma=0,08,~H=0.85.$}
\label{fig:label}
\end{figure}

\begin{table}[ht]
\centering
\caption{the mean value and standard deviation of the estimator $\mu$  }
\begin{tabular}{c|c|c|c|c}
\hline
 $\mu_0$ &0.6 & 1 & 1.5 &1.75 \\ \hline
 $H=0.65$\\ \hline
Mean value &0.6117876& 1.027482& 1.499406& 1.743383 \\ \hline
standard deviation &0.3219629& 0.3546342& 0.3302345&0.3440518\\
\hline
 $H=0.75$\\
\hline
Mean value &0.5902069& 0.9973574& 1.499907& 1.747189 \\ \hline
standard deviation &0.1488002& 0.1519977& 0.1513625&0.1489271\\
\hline
 $H=0.85$\\
\hline
Mean value &0.603559& 0.9996678& 1.498055& 1.750468 \\ \hline
standard deviation &0.061659& 0.0598438& 0.05994465&0.06095517\\
\hline
\end{tabular}
\end{table}

\begin{table}[ht]
\centering
\caption{the mean value and standard deviation of the estimator $\theta$  }
\begin{tabular}{c|c|c|c|c}
\hline
 $\theta_0$ &-0.7 & -0.8 & -0.9 &-0.95 \\ \hline
 $H=0.55$\\ \hline
Mean value &-0.69763& -0.795327& -0.9010633& -0.9515562 \\ \hline
standard deviation &0.03804792& 0.03779723& 0.03636855&0.03724209\\
\hline
 $H=0.65$\\
\hline
Mean value &-0.6983051& -0.8003717& -0.9000533&-0.9491352 \\ \hline
standard deviation &0.0260658& 0.0269981& 0.02682298&0.02685726\\
\hline
 $H=0.75$\\
\hline
Mean value &-0.7010637& -0.7988964& -0.900149& -0.9486839 \\ \hline
standard deviation &0.01824956& 0.01887855& 0.01815953&0.01940864\\
\hline
\end{tabular}
\end{table}

\section*{\normalsize  Conclusion and future work }

\quad\quad In this article, we mainly focus on two kinds of LSE of Vasicek-type stochastic differential equation driven by sub-fBm from discrete observations. In theorems 3.1 and 4.2, the consistency and asymptotic distribution are established. According to the research basis of this paper, we can also study other properties of the LSE of the Vasicek-type stochastic differential equations with discrete observations in the future.

\end{document}